\documentclass{amsart}

\usepackage{amssymb,amsmath,amsthm,latexsym,booktabs}

\usepackage{bbding}

\theoremstyle{definition}
\newtheorem{definition}{Definition}

\newtheorem{example}[definition]{Example}

\newtheorem{remark}[definition]{Remark}

\theoremstyle{plain}
\newtheorem{lemma}[definition]{Lemma}

\newtheorem{theorem}[definition]{Theorem}

\renewcommand{\lq}{\{}
\renewcommand{\rq}{\}}
\newcommand{\lqq}{\{}
\newcommand{\rqq}{\}'}

\newcommand{\aki}{\mathcal{A}}
\newcommand{\qua}{\mathcal{Q}}
\newcommand{\tpa}{\mathcal{T}}
\newcommand{\fpa}{\mathcal{F}}

\begin{document}

\title[Tangent algebras of monoassociative loops]
{Polynomial identities for tangent algebras of monoassociative loops}

\author[Bremner]{Murray R. Bremner}

\address{Department of Mathematics and Statistics,
University of Saskatchewan, Canada}

\email{bremner@math.usask.ca}

\author[Madariaga]{Sara Madariaga}

\address{Departamento de Matem\'aticas y Computaci\'on, Universidad de La Rioja,
Espa\~na}

\email{sara.madariaga@unirioja.es}

\begin{abstract}
We introduce degree $n$ Sabinin algebras, which are defined by the polynomial identities
up to degree $n$ in a Sabinin algebra.
Degree 4 Sabinin algebras can be characterized by the polynomial identities satisfied by
the commutator, associator and two quaternators in the free nonassociative algebra.
We consider these operations in a free power associative algebra and show that one
of the quaternators is redundant.
The resulting algebras provide the natural structure on the tangent space at the identity
element of an analytic loop for which all local loops satisfy monoassociativity,
$a^2 a \equiv a a^2$.
These algebras are the next step beyond Lie, Malcev, and Bol algebras.
We also present an identity of degree 5 which is satisfied by these three operations but
which is not implied by the identities of lower degree.
\end{abstract}

\maketitle

%%%%%%%%%%%%%%%%%%%%%%%%%%%%%%%%%%%%%%%%%%%%%%%%%%%%%%%%%%%%%%%%%%%%%%%%%%%%%%%%%%%%

\section{Introduction}

A Lie algebra is the natural algebraic structure on the tangent space at the
identity element of a Lie group.
In the theory of Lie groups, the third (or converse) Lie theorem is of great importance:
it states that a Lie algebra uniquely determines a local Lie group up to isomorphism.
If we remove the associativity condition on the multiplication in a Lie group,
we obtain the notion of an analytic loop.
Special classes of analytic loops are Moufang loops, Bol loops, and monoassociative loops.
The weakened associativity conditions defining these varieties of loops are equivalent to certain closure conditions
on the corresponding three-webs as described in Table \ref{looptable};
see \cite{AkivisGoldberg2,AkivisGoldberg3,AkivisShelekhovBook,Chern,Nagy}.

\begin{table}[h]
\begin{center}
\begin{tabular}{lll}
\toprule
LOCAL LOOPS & LOOP IDENTITIES & WEB CLOSURE \\
\midrule
abelian Lie groups & $u \cdot v = v \cdot u, (u \cdot v)\cdot w = u\cdot (v \cdot w)$ & (T) Thomsen \\
\midrule
Lie groups & $ (u \cdot v)\cdot w = u\cdot (v \cdot w) $ & (R) Reidemeister \\
\midrule
Moufang loops & $ (u \cdot v)\cdot (w\cdot u) = u\cdot ((v \cdot w)\cdot u) $ & (M) = (B$_l$) $\cap$ (B$_r$) \\
\midrule
left Bol loops & $(u\cdot(v\cdot u))\cdot w = u\cdot(v\cdot(u\cdot w))$ & (B$_l$) left Bol \\
\midrule
right Bol loops & $w \cdot ( ( u \cdot v ) \cdot u ) ) = ( ( w\cdot u ) \cdot v ) \cdot u$ & (B$_r$) right Bol \\
\midrule
monoassociative loops & $u^2 \cdot u = u \cdot u^2$ & (H) hexagonal \\
\bottomrule
\end{tabular}
\end{center}
\medskip
\caption{Varieties of loops and closure conditions}
\label{looptable}
\end{table}

Malcev (respectively Bol) algebras arise as tangent algebras of analytic Moufang (respectively Bol) loops.
A generalization of Lie's third theorem holds in these cases:
a Malcev (Bol) algebra determines a local Moufang (Bol) loop uniquely up to isomorphism.
In addition to a binary operation,
Bol algebras require a ternary operation in order to satisfy the generalization of Lie's third theorem.
Lie, Malcev, and Bol algebras can also be defined by the polynomial
identities of low degree satisfied by the commutator and associator in free
associative, alternative, and right alternative algebras; see \cite{KS,PIBol,PIS}.

In a free nonassociative algebra, the commutator and associator satisfy the Akivis identity;
this binary-ternary structure is called an Akivis algebra.
In general, for analytic loops, an Akivis algebra does not determine a local loop uniquely.
For example, for monoassociative loops, the prolonged Akivis algebra
defined in a fourth-order neighbourhood determines the local loop uniquely; see  \cite{AkivisGoldberg2}.
This prolonged Akivis algebra has two additional quaternary operations, called quaternators.
Mikheev \cite{Mikheev} proved that for monoassociative loops, one of the quaternators can be
expressed in terms of the other and the binary and ternary operations, because in this case
the local algebras are third-power associative.

In general, one considers an arbitrary analytic loop and the multilinear operations induced
by its multiplication on the tangent space at the identity element.
Akivis \cite{Akivis2} introduced the notion of closed $G$-structures on a differentiable quasigroup,
defined by a finite number of structure constants.
A local Akivis algebra, or its prolonged algebra, uniquely defines a differentiable quasigroup
if and only if the corresponding $G$-structure is closed.
The additional operations of higher degree in the prolonged algebra come from the higher-order
terms in the Taylor decomposition of the multiplication in the quasigroup.

If we include all these higher-degree operations, we obtain the structure described geometrically
in \cite{MikheevSabinin}.
The algebraic properties of these objects were clarified in \cite{PerezIzquierdo,ShestakovUmirbaev},
leading to the definition of Sabinin algebras, with an infinite family of operations which are
primitive elements with respect to the natural coproduct on the free nonassociative algebra.

Recall that a loop is Moufang if and only if all of its isotopes are alternative;
a loop is left (right) Bol if and only if all of its isotopes are left (right) alternative.
If a loop and all of its isotopes are third-power associative, then the loop is power-associative,
but not conversely: there are power-associative loops with isotopes which are not third-power associative \cite[Theorem 3.7]{Aczel}.
The corresponding tangent algebras can be defined by the polynomial identities of low degree
satisfied by certain multilinear operations in an underlying nonassociative algebra; see Table \ref{algebratable}.

\begin{table}[h]
\begin{center}
\begin{tabular}{ccc}
\toprule
UNDERLYING & PRIMITIVE & TANGENT  \\
ALGEBRA & OPERATIONS & ALGEBRA \\
\midrule
free associative & commutator & Lie algebras \\
\midrule
free alternative & \begin{tabular}{c} commutator, \\ (associator) \end{tabular} & Malcev algebras \\
\midrule
free right alternative & \begin{tabular}{c} commutator, \\ associator \end{tabular} & Bol algebras \\
\midrule
free power associative & \begin{tabular}{c} commutator, \\ associator, \\ quaternator \end{tabular} & BTQ algebras \\
\bottomrule
\end{tabular}
\end{center}
\medskip
\caption{Free algebras, multilinear operations, tangent algebras}
\label{algebratable}
\end{table}

The goal of the present paper is to determine the polynomial identities defining tangent algebras of monoassociative loops,
which we call BTQ algebras.
These structures have previously been investigated from a purely geometric point of view in
\cite{AkivisGoldberg2,AkivisShelekhov,Mikheev,Shelekhov}.
We will study the polynomial identities satisfied by these structures from the point of view of computer algebra.

We summarize the theoretical aspects of our computational methods.
Let $F\{X\}$ be the free nonassociative algebra on the set $X = \{ x_1, \dots, x_n  \}$ of generators over the field $F$,
and let $B_n$ be the subspace of $F\{X\}$ consisting of the multilinear polynomials of degree $n$.
Since every multilinear monomial consists of a permutation of $x_1 \cdots x_n$ together with an association type
(placement of parentheses), we have
  \[
  \dim B_n = n! K_n, \qquad
  K_n = \frac1n \binom{2n-2}{n-1} \quad \text{(Catalan number)}.
  \]
Let $F\{ \Omega; X \}$ be the free multioperator algebra on $X$ over $F$ where $\Omega$ is a finite set of
multilinear operations, and let $A_n$ be the subspace of $F\{ \Omega; X \}$ consisting of the multilinear elements
of degree $n$.
If $\omega \in \Omega$ is a $d$-ary operation, then we identify $\omega$ with an element of $p_\omega \in B_d$.
For example, in this paper an anticommutative binary operation corresponds to the commutator $[a,b] = ab - ba$,
and a ternary operation corresponds to the associator $(a,b,c) = (ab)c - a(bc)$; we will also consider two
quaternary operations (see Definition \ref{defquaternators}).
For each $n$ we consider a linear map $E_n \colon A_n \to B_n$, which we call the expansion map, defined
by replacing each occurrence of $\omega$ by $p_\omega$ and making the appropriate substitutions.  For example,
  \[
  E_4( [ (a,b,c), d ] )
  =
  ((ab)c)d - (a(bc))d - d((ab)c) + d(a(bc)).
  \]
The basic computational principle is that the kernel of $E_n$ consists of the polynomial identities of degree $n$
satisfied by the multilinear operations defined by the elements $\{ p_\omega \mid \omega \in \Omega \}$ in $F\{X\}$.

We modify this construction as follows.  We replace the free nonassociative algebra $F\{X\}$
by the free algebra $F_\mathcal{V}\{X\}$ in some variety $\mathcal{V}$ of nonassociative algebras;
for example, associative, alternative, right alternative, or power associative.
Then $B_{n,\mathcal{V}}$ is the multilinear subspace of degree $n$ in $F_\mathcal{V}\{X\}$,
which we can identify with the quotient space of $B_n$ by the subspace of multilinear identities of
degree $n$ satisfied by the variety $\mathcal{V}$.
The kernel of $E_{n,\mathcal{V}} \colon A_n \to B_{n,\mathcal{V}}$ consists of the polynomial identities of degree $n$
satisfied by the elements $\{ p_\omega \mid \omega \in \Omega \}$ in $F_\mathcal{V}\{X\}$.
Referring to Table \ref{algebratable}, we see that $F_\mathcal{V}\{X\}$ is the underlying algebra,
the elements $p_\omega$ for $\omega \in \Omega$ are the primitive operations, and
the polynomial identities in the kernel of $E_{n,\mathcal{V}}$ define the variety of tangent algebras.

%%%%%%%%%%%%%%%%%%%%%%%%%%%%%%%%%%%%%%%%%%%%%%%%%%%%%%%%%%%%%%%%%%%%%%%%%%%%%%%%%%%%

\section{Quaternary operations and Sabinin algebras}

We consider the free nonassociative algebra $F\{X\}$ generated by the set $X$ over
the field $F$. We define a coproduct $\Delta\colon F\{X\} \to F\{X\} \otimes F\{X\}$
inductively on the degree of basis monomials in $F\{X\}$ by setting $\Delta(x) = x
\otimes 1 + 1 \otimes x$ for $x \in X$ and assuming that $\Delta$ is an algebra
morphism: $\Delta(fg) = \Delta(f) \Delta(g)$.  Multilinear operations
correspond to nonassociative polynomials $f \in F\{X\}$ which are primitive with respect 
to this coproduct:
  \[
  \Delta(f ) = f \otimes 1 + 1 \otimes f.
  \]
Shestakov and Umirbaev \cite{ShestakovUmirbaev} give an algorithm for constructing a complete
set of these primitive elements. 

\begin{definition}
Let $f, g_1, \dots, g_k \in F\{X\}$ be multilinear elements of the same degree.
We say that $f$ is a \textbf{consequence} of $g_1, \dots, g_k$ if it is a linear
combination of the elements obtained by permuting the variables in $g_1, \dots, g_k$.
\end{definition}

\begin{example}
In degree 2 every primitive element is a consequence of the commutator $[a,b] = ab - ba$.
In degree 3 every primitive element is a consequence of the iterated commutator $[[a,b],c]$ and the associator $(a,b,c) = (ab)c - a(bc)$.
\end{example}

\begin{definition}
By an \textbf{Akivis element} in the free nonassociative algebra, we mean a polynomial
which can be expressed using only the commutator and the associator.
\end{definition}

It is easy to verify that every Akivis element is primitive.
In degree 4, every multilinear Akivis element is a consequence of these six polynomials:
  \begin{equation}
  \label{akivis}
  \begin{array}{ll}
  \aki_1(a,b,c,d) = [ [ [ a, b ], c ], d ],
  &\qquad
  \aki_2(a,b,c,d) = [ ( a, b, c ), d ],
  \\[4pt]
  \aki_3(a,b,c,d) = [ [ a, b ], [ c, d ] ],
  &\qquad
  \aki_4(a,b,c,d) = ( [ a, b ], c, d ),
  \\[4pt]
  \aki_5(a,b,c,d) = ( a, [ b, c ], d ),
  &\qquad
  \aki_6(a,b,c,d) = ( a, b, [ c, d ] ).
  \end{array}
  \end{equation}
In degree 4 there are two primitive elements which are not Akivis elements:
  \begin{equation}
  \label{pijformula}
  \begin{array}{l}
  p_{2,1}( a b, c, d )
  =
  ( a b, c, d ) - a ( b, c, d ) - b ( a, c, d ),
  \\[4pt]
  p_{1,2}( a, b c, d )
  =
  ( a, b c, d ) - b ( a, c, d ) - c ( a, b, d ).
  \end{array}
  \end{equation}
Every primitive multilinear nonassociative polynomial of degree 4 is a consequence of the previous eight elements.
We will use two slightly modified quaternary operations,
which can be regarded as measuring the deviation of the associators $(-,c,d)$ and $(a,-,d)$ from being derivations.

\begin{definition} \label{defquaternators}
The two \textbf{quaternators} are defined as follows:
  \begin{equation}
  \label{quaternators}
  \begin{array}{l}
  \qua_1(a,b,c,d) = ( ab, c, d ) - a ( b, c, d ) - ( a, c, d ) b,
  \\[4pt]
  \qua_2(a,b,c,d) = ( a, bc, d ) - b ( a, c, d ) - ( a, b, d ) c,
  \end{array}
  \end{equation}
It is clear that $\qua_1$ and $\qua_2$ are equivalent to $p_{2,1}$ and $p_{1,2}$ modulo $\aki_2$.
\end{definition}

We now recall some basic definitions and results from \cite{ShestakovUmirbaev}.
All nonassociative monomials are assumed to be right-normed:
$x_1 x_2 x_3 \cdots x_m = ( \cdots ( ( x_1 x_2 ) x_3 ) \cdots x_m )$.

\begin{definition}
The ordered pair $(u_1,u_2)$ is a \textbf{2-decomposition} of $u = x_1 \cdots
x_m$ if $u_1 = x_{i_1} \cdots x_{i_k}$ and $u_2 = x_{i_{k+1}} \cdots x_{i_m}$
where the subsets $I_1 = \{ i_1, \dots, i_k \}$ and $I_2 = \{ i_{k+1}, \dots,
i_m \}$ satisfy
  \[
  i_1 < \cdots < i_k,
  \qquad
  i_{k+1} < \cdots < i_m,
  \qquad
  I_1 \cup I_2 = \{ 1, \dots, m \},
  \qquad
  I_1 \cap I_2 = \emptyset.
  \]
If $k = 0$ then $u_1 = 1$; if $k = m$ then $u_2 = 1$. (More generally, see \cite[p.~538]{ShestakovUmirbaev}.)
\end{definition}

\begin{definition} \cite[p.~544]{ShestakovUmirbaev}
We inductively define the following polynomials in the free nonassociative algebra; these elements are primitive:
  \allowdisplaybreaks
  \begin{align*}
  &
  p_{1,1}( x, y, z )
  =
  ( x, y, z ), \quad \text{(associator)}
  \\
  &
  p_{m,n}( x_1 \cdots x_m, y_1 \cdots y_n, z )
  \\
  &=
  ( \, x_1 \cdots x_m, \, y_1 \cdots y_n, \, z \, )
  -
  \!\!\!\!\!\!\!\!\!\!\!\!
  \sum_{\tiny
  \begin{array}{c}
  u_1 \ne 1 \; \text{or} \; v_1 \ne 1 \\
  u_2 \ne 1 \; \text{and} \; v_2 \ne 1
  \end{array}}
  \!\!\!\!\!\!\!\!\!\!\!\!
  u_1 v_1 p_{m-\deg(u_1),n-\deg(v_1)}
  ( u_2, v_2, z ),
  \end{align*}
where $m + n \ge 3$ and the sum is over all 2-decompositions of $u = x_1 \cdots
x_m$ and $v = y_1 \cdots y_n$.
\end{definition}

\begin{remark}
We compute the polynomials of degree 4. For $m = 2$, $n = 1$ we have $u = x_1
x_2$, $v = y$ and the 2-decompositions are
  \[
  ( u_1, u_2 ) = ( x_1 x_2, 1 ), ( x_1, x_2 ), ( x_2, x_1 ), ( 1, x_1 x_2 ),
  \qquad
  ( v_1, v_2 ) = ( y, 1 ), ( 1, y ).
  \]
We may only combine $( x_1, x_2 )$ with $( 1, y )$, and $( x_2,
x_1 )$ with $( 1, y )$. We obtain
  \allowdisplaybreaks
  \begin{align*}
  &
  p_{2,1}( x_1 x_2, y, z )
  =
  ( x_1 x_2, y, z ) - x_1 p_{1,1}( x_2, y, z ) - x_2 p_{1,1}( x_1, y, z )
  \\
  &=
  ( ( x_1 x_2 ) y ) z - ( x_1 x_2 ) ( y z )
  - x_1 ( ( x_2 y ) z ) + x_1 ( x_2 ( y z ) )
  - x_2 ( ( x_1 y ) z ) + x_2 ( x_1 ( y z ) ).
  \end{align*}
Similarly, for $m = 1$, $n = 2$ we obtain
  \allowdisplaybreaks
  \begin{align*}
  &
  p_{1,2}( x, y_1 y_2, z )
  =
  ( x, y_1 y_2, z ) - y_1 p_{1,1}( x, y_2, z ) - y_2 p_{1,1}( x, y_1, z )
  \\
  &=
  ( x ( y_1 y_2 ) ) z - x ( ( y_1 y_2 ) z )
  - y_1 ( ( x y_2 ) z ) + y_1 ( x ( y_2 z ) )
  - y_2 ( ( x y_1 ) z ) + y_2 ( x ( y_1 z ) ).
  \end{align*}
\end{remark}

\begin{definition}
\cite[p.~539]{ShestakovUmirbaev}
In any nonassociative algebra, we define the multilinear \textbf{Sabinin operations}
in terms of the commutator and the polynomials $p_{m,n}$:
  \allowdisplaybreaks
  \begin{align*}
  \langle y, z \rangle &= - [ y, z ],
  \\
  \langle x_1, \dots, x_m, y, z \rangle
  &=
  - p_{m,1}( u, y, z ) + p_{m,1}( u, z, y ),
  \end{align*}
where $u = x_1 \cdots x_m$ and $m \ge 1$, and
  \allowdisplaybreaks
  \begin{align*}
  &
  \Phi_{m,n}( x_1, \dots, x_m; y_1, \dots, y_n)
  =
  \\
  &
  \frac{1}{m!n!}
  \sum_{\sigma \in S_m}
  \sum_{\tau \in S_n}
  p_{m,n-1}( x_{\sigma(1)}, \dots, x_{\sigma(m)}; y_{\tau(1)}, \dots, y_{\tau(n-1)};y _{\tau(n)} ),
  \end{align*}
where $m \ge 1$ and $n \ge 2$.
\end{definition}

\begin{remark}
In degrees 2, 3 and 4 we have the following operations:
  \allowdisplaybreaks
  \begin{align*}
  \langle y, z \rangle &= - [ y, z ],
  \\
  \langle x, y, z \rangle
  &=
  - p_{1,1}( x, y, z ) + p_{1,1}( x, z, y ),
  \\
  \Phi_{1,2}( x,y, z )
  &=
  \tfrac12 \big( p_{1,1}( x, y, z ) + p_{1,1}( x, z, y ) \big),
  \\
  \langle w, x, y, z \rangle
  &=
  - p_{2,1}( wx, y, z ) + p_{2,1}( wx, z, y ),
  \\
  \Phi_{1,3}( w, x, y, z )
  &=
  \tfrac16
  \big(
  p_{1,2}( w, xy, z )
  +
  p_{1,2}( w, xz, y )
  +
  p_{1,2}( w, yx, z )
  \\
  &\quad\quad\quad
  +
  p_{1,2}( w, yz, x )
  +
  p_{1,2}( w, zx, y )
  +
  p_{1,2}( w, zy, x )
  \big),
  \\
  \Phi_{2,2}( w, x, y, z )
  &=
  \tfrac14
  \big(
  p_{2,1}( wx, y, z )
  +
  p_{2,1}( wx, z, y )
  \\
  &\quad\quad\quad
  +
  p_{2,1}( xw, y, z )
  +
  p_{2,1}( xw, z, y )
  \big).
  \end{align*}
Using equations \eqref{pijformula}, and sorting the terms first by association type in the order
  \[
  ((ab)c)d, \quad
  (a(bc))d, \quad
  (ab)(cd), \quad
  a((bc)d), \quad
  a(b(cd)),
  \]
and then by permutation of the letters, we obtain the fully expanded quaternary operations:
  \allowdisplaybreaks
  \begin{align*}
  &\langle w, x, y, z \rangle =
  \\
  &
  - ( ( w x ) y ) z
  + ( ( w x ) z ) y
  + ( w x ) ( y z )
  - ( w x ) ( z y )
  + w ( ( x y ) z )
  - w ( ( x z ) y )
  \\
  &
  + x ( ( w y ) z )
  - x ( ( w z ) y )
  - w ( x ( y z ) )
  + w ( x ( z y ) )
  - x ( w ( y z ) )
  + x ( w ( z y ) ),
  \\[6pt]
  &
  \Phi_{1,3}( w, x, y, z ) =
  \\
  &
  \tfrac16
  \big(
  ( w ( x y ) ) z
  + ( w ( x z ) ) y
  + ( w ( y x ) ) z
  + ( w ( y z ) ) x
  + ( w ( z x ) ) y
  + ( w ( z y ) ) x
  \\
  &
  - w ( ( x y ) z )
  - w ( ( x z ) y )
  - w ( ( y x ) z )
  - w ( ( y z ) x )
  - w ( ( z x ) y )
  - w ( ( z y ) x )
  \\
  &
  - 2 x ( ( w y ) z )
  - 2 x ( ( w z ) y )
  - 2 y ( ( w x ) z )
  - 2 y ( ( w z ) x )
  - 2 z ( ( w x ) y )
  - 2 z ( ( w y ) x )
  \\
  &
  + 2 x ( w ( y z ) )
  + 2 x ( w ( z y ) )
  + 2 y ( w ( x z ) )
  + 2 y ( w ( z x ) )
  + 2 z ( w ( x y ) )
  + 2 z ( w ( y x ) )
  \big),
  \\[6pt]
  &
  \Phi_{2,2}( w, x, y, z ) =
  \\
  &
  \tfrac14
  \big(
    ( ( w x ) y ) z
  + ( ( w x ) z ) y
  + ( ( x w ) y ) z
  + ( ( x w ) z ) y
  - ( w x ) ( y z )
  - ( w x ) ( z y )
  \\
  &
  - ( x w ) ( y z )
  - ( x w ) ( z y )
  - 2 w ( ( x y ) z )
  - 2 w ( ( x z ) y )
  - 2 x ( ( w y ) z )
  - 2 x ( ( w z ) y )
  \\
  &
  + 2 w ( x ( y z ) )
  + 2 w ( x ( z y ) )
  + 2 x ( w ( y z ) )
  + 2 x ( w ( z y ) )
  \big).
  \end{align*}
These three Sabinin operations are equivalent to the two quaternators of equation \eqref{quaternators}
in the following sense.  Consider the multilinear subspace $N$ of degree 4 in the
free nonassociative algebra on 4 generators, the subspace $P$ of primitive
elements, and the subspace $A$ of Akivis elements.  Clearly $A \subseteq P
\subseteq N$, and all these spaces have a natural action of the symmetric group
$S_4$ by permutation of the variables.  We have two sets of generators for the
quotient module $P/A$: the cosets of the two quaternators, and the cosets of
the three Sabinin operations.
\end{remark}

\begin{lemma}
The Sabinin operations in degree 4 can be expressed in terms of the Akivis
elements and the quaternators as follows:
  \allowdisplaybreaks
  \begin{align*}
  \langle a, b, c, d \rangle
  &=
  - \aki_2( a, c, d, b )
  + \aki_2( a, d, c, b )
  - \qua_1( a, b, c, d )
  + \qua_1( a, b, d, c ),
  \\[4pt]
  \Phi_{1,3}( a, b, c, d )
  &=
   2 \aki_2( a, b, c, d )
  +2 \aki_2( a, b, d, c )
  +2 \aki_2( a, c, b, d )
  -4 \aki_2( a, c, d, b )
  \\
  &\quad
  +2 \aki_2( a, d, b, c )
  +2 \aki_2( a, d, c, b )
  -3 \aki_5( a, b, c, d )
  -  \aki_5( a, b, d, c )
  \\
  &\quad
  -  \aki_5( a, c, d, b )
  -2 \aki_6( a, b, c, d )
  -2 \aki_6( a, c, b, d )
  -2 \qua_1( a, b, c, d )
  \\
  &\quad
  +2 \qua_1( a, b, d, c )
  -2 \qua_1( a, c, b, d )
  +2 \qua_1( a, c, d, b )
  +6 \qua_2( a, b, c, d ),
  \\[4pt]
  \Phi_{2,2}( a, b, c, d )
  &=
   2 \aki_2( a, c, d, b )
  +2 \aki_2( a, d, c, b )
  -  \aki_4( a, b, c, d )
  -  \aki_4( a, b, d, c )
  \\
  &\quad
  +2 \qua_1( a, b, c, d )
  +2 \qua_1( a, b, d, c ).
  \end{align*}
The quaternators can be expressed in terms of the Akivis elements and the
Sabinin operations in degree 4 as follows:
  \allowdisplaybreaks
  \begin{align*}
  \qua_1(a,b,c,d)
  =
  & - \aki_2(a,c,d,b)
  + \tfrac{1}{4} \aki_4(a,b,c,d)
  + \tfrac{1}{4} \aki_4(a,b,d,c)
  \\
  &
  - \tfrac{1}{2} \langle a,b,c,d \rangle
  + \tfrac{1}{4} \Phi_{2,2}(a,b,c,d)
  \\[4pt]
  \qua_2(a,b,c,d)
  =
  & -\tfrac{1}{3} \aki_2(a,b,c,d)
  - \tfrac{2}{3} \aki_2(a,b,d,c)
  -\tfrac{1}{3} \aki_2(a,c,b,d)
  \\
  &
  + \tfrac{1}{3} \aki_2(a,c,d,b)
  + \tfrac{1}{2} \aki_5(a,b,c,d)
  + \tfrac{1}{6} \aki_5(a,b,d,c) \\
  &
  + \tfrac{1}{6} \aki_5(a,c,d,b)
  + \tfrac{1}{3} \aki_6(a,b,c,d)
  + \tfrac{1}{3} \aki_6(a,c,b,d)
  \\
  &
  -\tfrac{1}{3} \langle a,b,c,d \rangle
  - \tfrac{1}{3} \langle a,c,b,d \rangle
  + \tfrac{1}{6} \Phi_{1,3}(a,b,c,d)
  \end{align*}
\end{lemma}

\begin{proof}
Expand both sides of each equation in the free nonassociative algebra.
\end{proof}

\begin{definition} \label{definitionsabinin}
\cite[p.~544]{ShestakovUmirbaev}
A \textbf{Sabinin algebra} is a vector space with two infinite families of multilinear operations
  \begin{align*}
  &
  \left< x_1, x_2, \dots, x_m; y, z \right>, \quad m \ge 0,
  \\
  &
  \Phi( x_1, x_2, \dots, x_m; y_1, y_2, \dots, y_n ), \quad m \ge 1, n \ge 2,
  \end{align*}
satisfying the polynomial identities
  \begin{align}
  &
  \left< x_1, x_2, \dots, x_m; y, z \right>
  +
  \left< x_1,x_2,\dots,x_m;z,y \right>
  \equiv
  0,
  \label{S1}
  \\
  &
  \left< x_1, x_2, \dots, x_r, a, b, x_{r+1}, \dots, x_m; y, z \right>
  \label{S2}
  \\
  &\quad
  -
  \left< x_1, x_2, \dots, x_r, b, a, x_{r+1}, \dots, x_m; y, z \right>
  \notag
  \\
  &\quad
  +
  \sum_{k=0}^r
  \sum_{\alpha}
  \left< x_{\alpha_1}, \dots, x_{\alpha_k},
  \left< x_{\alpha_{k+1}}, \dots, x_{\alpha_r}; a, b \right>,
  \dots,
  x_m; y, z \right>
  \equiv
  0,
  \notag
  \\
  &
  \sigma_{x,y,z}
  \Big(
  \left< x_1, \dots, x_r, x; y, z \right>
  \label{S3}
  \\
  &\quad
  +
  \sum_{k=0}^r
  \sum_{\alpha}
  \left< x_{\alpha_1}, \dots, x_{\alpha_k};
  \left< x_{\alpha_{k+1}}, \dots, x_{\alpha_r}; y, z \right>, x
  \right>
  \Big)
  \equiv
  0,
  \notag
  \\
  &
  \Phi( x_1, \dots, x_m; y_1, \dots, y_n )
  \equiv
  \Phi( x_{\tau(1)}, \dots, x_{\tau(m)}; y_{\delta(1)}, \dots, y_{\delta(n)} ),
  \label{S4}
  \end{align}
where $\alpha \in S_r$ $(i \mapsto \alpha_i)$, $\alpha_1 < \alpha_2 < \cdots
<\alpha_k$, $\alpha_{k+1} < \cdots < \alpha_r$, $k=0,1,\dots,r$, $r \geq 0$,
$\sigma_{x,y,z}$ denotes the cyclic sum over $x,y,z$, and $\tau \in S_m$,
$\delta \in S_n$.
\end{definition}

We introduce the following varieties of multioperator algebras.

\begin{definition}
A \textbf{degree $d$ Sabinin algebra} for $d \ge 2$ is obtained from Definition \ref{definitionsabinin}
by considering only the multilinear operations of degree up to $d$, together with the polynomial
identities of degree up to $d$ involving these operations.
\end{definition}

\begin{remark}
For $d = 2$ we obtain the variety of anticommutative algebras
with one bilinear operation $[a,b]$ satisfying $[a,b] + [b,a] \equiv 0$.  For
$d = 3$ we obtain the variety of algebras with a bilinear operation $[a,b]$ and
two trilinear operations $\langle a,b,c \rangle$ and $\Phi_{1,2}(a,b,c)$
satisfying the identities
  \begin{align*}
  &
  [a,b] + [b,a] \equiv 0,
  \\
  &
  \langle a, b, c \rangle + \langle a, c, b \rangle \equiv 0,
  \\
  &
  \langle a, b, c \rangle + [[b,c],a] +
  \langle b, c, a \rangle + [[c,a],b] +
  \langle c, a, b \rangle + [[a,b],c]
  \equiv 0,
  \\
  &
  \Phi_{1,2}(a,b,c) \equiv \Phi_{1,2}(a,c,b).
  \end{align*}
It is easy to verify that the associator can be recovered by the formula
  \[
  (a,b,c) = \tfrac12 \big( 2 \Phi_{1,2}(a,b,c) - \langle a,b,c \rangle \big),
  \]
and so this variety is equivalent to the variety of Akivis algebras.
\end{remark}

%%%%%%%%%%%%%%%%%%%%%%%%%%%%%%%%%%%%%%%%%%%%%%%%%%%%%%%%%%%%%%%%%%%%%%%%%%%%%%%%%%%%

\section{Degree 4 Sabinin algebras (BTQQ algebras)}

In this section we consider degree 4 Sabinin algebras and find a simpler equivalent variety.
We first consider the identities relating the operations of degrees 2, 3, 4:
  \[
  \langle y, z \rangle,
  \quad
  \langle x, y, z \rangle,
  \quad
  \langle w, x, y, z \rangle,
  \quad
  \Phi_{1,2}( x, y, z ),
  \quad
  \Phi_{1,3}( w, x, y, z ),
  \quad
  \Phi_{2,2}( w, x, y, z ).
  \]
From identity \eqref{S1} we obtain skew-symmetry in the last two arguments of
the first three operations:
  \[
  \langle y, z \rangle + \langle z, y \rangle \equiv 0,
  \quad
  \langle x, y, z \rangle + \langle x, z, y \rangle \equiv 0,
  \quad
  \langle w, x, y, z \rangle + \langle w, x, z, y \rangle \equiv 0.
  \]
In identity \eqref{S2}, we must take $r = m = 0$, and we obtain
  \[
  \left< a, b, y, z \right>
  -
  \left< b, a, y, z \right>
  +
  \left< \left< a, b \right>, y, z \right>
  \equiv
  0.
  \]
In identity \eqref{S3}, we can take $r = 0$ or $r = 1$.  For $r = 0$, we obtain
the Akivis identity,
  \[
  \langle x, y, z \rangle
  +
  \langle y, z, x \rangle
  +
  \langle z, x, y \rangle
  +
  \langle \langle x, y \rangle, z \rangle
  +
  \langle \langle y, z \rangle, x \rangle
  +
  \langle \langle z, x \rangle, y \rangle
  \equiv
  0.
  \]
For $r = 1$, we obtain
  \[
  \sigma_{x,y,z}
  \big(
  \left< x_1, x, y, z \right>
  +
  \left< \left< x_1, y, z \right>, x \right>
  +
  \left< x_1, \left< y, z \right>, x \right>
  \big)
  \equiv
  0.
  \]
From identity \eqref{S4}, we obtain symmetries of the last three operations:
  \allowdisplaybreaks
  \begin{align*}
  &
  \Phi_{1,2}( x, y, z )
  \equiv
  \Phi_{1,2}( x, z, y ),
  \\
  &
  \Phi_{1,3}( w, x, y, z )
  \equiv
  \Phi_{1,3}( w, x, z, y )
  \equiv
  \Phi_{1,3}( w, y, x, z )
  \equiv
  \Phi_{1,3}( w, y, z, x )
  \\
  &\qquad
  \equiv
  \Phi_{1,3}( w, z, x, y )
  \equiv
  \Phi_{1,3}( w, z, y, x ),
  \\
  &
  \Phi_{2,2}( w, x, y, z )
  \equiv
  \Phi_{2,2}( w, x, z, y )
  \equiv
  \Phi_{2,2}( x, w, y, z )
  \equiv
  \Phi_{2,2}( x, w, z, y ).
  \end{align*}
The next result gives equivalent identities relating the commutator, associator, and the two quaternators.
We regard the quaternators as multilinear operations, and use the notation $\lq a,b,c,d \rq = \qua_1(a,b,c,d)$
and $\lqq a,b,c,d \rqq = \qua_2(a,b,c,d)$.

\begin{theorem} \label{BTQQtheorem}
Every polynomial identity of degree at most 4 satisfied by the commutator
$[a,b]$, the associator $(a,b,c)$, and the two quaternators $\lq a,b,c,d \rq$
and $\lqq a,b,c,d \rqq$, is a consequence of the following five independent
identities:
  \begin{align*}
  [a,b] + [b,a]
  &\equiv
  0,
  \\[4pt]
  [[a,b],c] + [[b,c],a] + [[c,a],b]
  &\equiv
  (a,b,c) - (a,c,b) - (b,a,c)
  \\
  &\quad
  + (b,c,a) + (c,a,b) - (c,b,a),
  \\[4pt]
  ( [a,b], c, d ) - [ a, (b,c,d) ] + [ b, (a,c,d) ]
  &\equiv
  \lq a,b,c,d \rq - \lq b,a,c,d \rq,
  \\[4pt]
  (a, [b, c], d ) - [b, (a, c, d) ] + [c, (a, b, d) ]
  &\equiv
  \lqq a, b, c, d \rqq - \lqq a, c, b, d \rqq,
  \\[4pt]
  [ b, (a,c,d) ] - [ b, (a,d,c) ] - ( a, b, [c,d] )
  &\equiv
  \lq a,b,c,d \rq
  - \lq a,b,d,c \rq
  \\
  &\quad
  - \lqq a,b,c,d \rqq
  + \lqq a,b,d,c \rqq.
  \end{align*}
\end{theorem}

\begin{proof}
We proved this result using the computer algebra system Maple.

(1) We consider the association types in degree 4 for these four operations:
  \begin{align*}
  &[ -, [ -, [ -, - ] ] ], &\;
  &[ -, ( -, -, - ) ], &\;
  &[ [ -, - ], [ -, - ] ], &\;
  &( -, -, [ -, - ] ),
  \\
  &( -, [ -, - ], - ), &\;
  &( [ -, - ], -, - ), &\;
  &\lq -, -, -, - \rq, &\;
  &\lqq -, -, -, - \rqq.
  \end{align*}
We call these the BTQQ types (binary, ternary, and two quaternary operations).
We generate the multilinear BTQQ monomials by applying all permutations of $a,b,c,d$ to the BTQQ types.
The skew-symmetry of the commutator implies that many of these monomials are scalar multiples
of other monomials; 
the space spanned by these monomials has a basis of $12 + 24 + 3 + 12 + 12 + 12 + 24 + 24 = 123$ distinct monomials.
We construct a matrix with a $123 \times 120$ left block $X$ and the identity matrix as the right block:
  \[
  \left[ \begin{array}{c|c} X & I \end{array} \right].
  \]
The rows of $X$ contain the coefficient vectors of the expansions of the multilinear BTQQ monomials into the free nonassociative algebra; 
the columns of $X$ correspond to the 24 multilinear monomials in each of the 5 association types in the free nonassociative algebra.
The columns of the right block correspond to the multilinear BTQQ monomials.
Each row of this matrix represents an equation saying that a BTQQ monomial equals its expansion in the free nonassociative algebra.
We compute the RCF (row canonical
form) of this matrix, and identify the rows whose leading 1s are in the right
block. There are 45 such rows, and they form a basis of the subspace of all polynomial identities
in degree 4 satisfied by the four operations, including the consequences of the
Akivis identity.  We sort these identities by increasing number of terms.

(2) We construct a matrix with a $123 \times 123$ upper block $A$ and a $24 \times 123$ lower block $B$,
initialized to zero:
  \[
  \left[ \begin{array}{c} A \\ \midrule B \end{array} \right].
  \]
We consider the consequences in degree 4 of the Akivis identity $A(a,b,c)$;
the last identity is redundant:
  \[
  A([a,d],b,c), \quad A(a,[b,d],c), \quad A(a,b,[c,d]), \quad [ A(a,b,c), d ], \quad [ d, A(a,b,c) ].
  \]
For each of these elements, we apply all permutations to the variables, store the coefficient
vectors in the rows of $B$, and compute the RCF; the lower block becomes zero after each iteration.
The result is a matrix of rank 10 whose row space is the subspace of all consequences in degree 4 of the Akivis identity.
Combining this with the result from step (1), we conclude that there is a quotient space of dimension $45 - 10 = 35$
consisting of the new identities in degree 4 satisfied by the four operations.

(3) Using the result of step (2), we take each basis vector from step (1),
apply all permutations of the variables, store the coefficient vectors in $B$,
and compute the RCF. Only four identities increase
the rank, and another computation shows that one of these is a consequence of the others
(and the liftings of the Akivis identity), so we have three new independent
identities:
    \begin{align*}
    &
    [c, [d, b, a]] - [d, [c, b, a]] - [[c, d], b, a] + \lq c, d, b, a \rq - \lq d, c, b, a \rq \equiv 0,
    \\[4pt]
    &
    [b, [c, d, a]] - [d, [c, b, a]] - [c, [b, d], a] + \lqq c, b, d, a \rqq - \lqq c, d, b, a \rqq \equiv 0,
    \\[4pt]
    &
    [c, [d, a, b]] - [d, [c, b, a]] - [d, c, [a, b]] - [[c, d], b, a]
    \\
    & \quad\quad\quad
    + \lq c, d, b, a \rq - \lq d, c, a, b\rq + \lqq d, c, a, b \rqq - \lqq d, c, b, a\rqq \equiv 0.
    \end{align*}
It is easy to verify that these identities are equivalent to the last three identities in the statement
of this Theorem.
\end{proof}

\begin{definition}
A \textbf{BTQQ algebra} is a vector space with a bilinear operation $[a,b]$, a
trilinear operation $(a,b,c)$, and two quadrilinear operations $\lq a,b,c,d
\rq$ and $\lqq a,b,c,d \rqq$, satisfying the identities of Theorem
\ref{BTQQtheorem}.
\end{definition}

\begin{theorem}
On every degree 4 Sabinin algebra, we can define operations which make it into
a BTQQ algebra, and conversely.
\end{theorem}

\begin{proof}
Given a BTQQ algebra, we define
\begin{align*}
  \langle a , b \rangle
  &=
  - [ a , b ],
  \\
  \langle a , b , c \rangle
  &=
  - ( a , b , c ) + ( a , c , b ),
  \\
  \Phi_{1,2}( a , b , c )
  &=
  ( a , b , c ) + ( a , c , b ),
  \\
  \langle a, b, c, d \rangle
  &=
  - \aki_2( a, c, d, b )
  + \aki_2( a, d, c, b )
  - \qua_1( a, b, c, d )
  + \qua_1( a, b, d, c ),
  \\
  \Phi_{1,3}( a, b, c, d )
  &=
   2 \aki_2( a, b, c, d )
  +2 \aki_2( a, b, d, c )
  +2 \aki_2( a, c, b, d )
  -4 \aki_2( a, c, d, b )
  \\
  &\quad
  +2 \aki_2( a, d, b, c )
  +2 \aki_2( a, d, c, b )
  -3 \aki_5( a, b, c, d )
  -  \aki_5( a, b, d, c )
  \\
  &\quad
  -  \aki_5( a, c, d, b )
  -2 \aki_6( a, b, c, d )
  -2 \aki_6( a, c, b, d )
  -2 \qua_1( a, b, c, d )
  \\
  &\quad
  +2 \qua_1( a, b, d, c )
  -2 \qua_1( a, c, b, d )
  +2 \qua_1( a, c, d, b )
  +6 \qua_2( a, b, c, d ),
  \\
  \Phi_{2,2}( a, b, c, d )
  &=
   2 \aki_2( a, c, d, b )
  +2 \aki_2( a, d, c, b )
  -  \aki_4( a, b, c, d )
  -  \aki_4( a, b, d, c )
  \\
  &\quad
  +2 \qua_1( a, b, c, d )
  +2 \qua_1( a, b, d, c ),
  \end{align*}
where the Akivis elements $\aki_i(a,b,c,d)$ are defined in equations \eqref{akivis}.
Straightforward computations prove that these operations satisfy all the identities of a degree 4 Sabinin algebra.
For example, for the Akivis identity, we have
  \begin{align*}
  &
  \langle x, y, z \rangle
  +
  \langle y, z, x \rangle
  +
  \langle z, x, y \rangle
  +
  \langle \langle x, y \rangle, z \rangle
  +
  \langle \langle y, z \rangle, x \rangle
  +
  \langle \langle z, x \rangle, y \rangle
  =
  \\
  &
  - (x,y,z) + (x,z,y) - (y,z,x) + (y,x,z) - (z,x,y) + (z,y,x) &
  \\
  &
  + [[x,y],z] + [[y,z],x] + [[z,x],y] = 0, &
\end{align*}
where the first equality uses the definitions of the operations $\langle -, - \rangle$ and $\langle -, -, - \rangle$
and the second equality follows from the second identity in Theorem \ref{BTQQtheorem}.

Conversely, given a degree 4 Sabinin algebra, we define
  \begin{align*}
  [ a , b ]
  &=
  - \langle a, b \rangle,
  \\
  ( a , b , c )
  &=
  \tfrac12 \big( 2 \Phi_{1,2} ( a , b , c ) - \langle a , b , c \rangle \big),
  \\
  \lq a,b,c,d \rq
  &=
  - \widetilde{\aki}_2(a,c,d,b) + \tfrac{1}{4} \widetilde{\aki}_4(a,b,c,d) + \tfrac{1}{4} \widetilde{\aki}_4(a,b,d,c)
  \\
  &\quad
  - \tfrac{1}{2} \langle a,b,c,d \rangle + \tfrac{1}{4} \Phi_{2,2}(a,b,c,d),
  \\
  \lqq a,b,c,d \rqq
  &=
  -\tfrac{1}{3} \widetilde{\aki}_2(a,b,c,d) - \tfrac{2}{3} \widetilde{\aki}_2(a,b,d,c) -\tfrac{1}{3} \widetilde{\aki}_2(a,c,b,d)
  \\
  &\quad
  + \tfrac{1}{3} \widetilde{\aki}_2(a,c,d,b) + \tfrac{1}{2} \widetilde{\aki}_5(a,b,c,d) + \tfrac{1}{6} \widetilde{\aki}_5(a,b,d,c)
  \\
  &\quad
  + \tfrac{1}{6} \widetilde{\aki}_5(a,c,d,b) + \tfrac{1}{3} \widetilde{\aki}_6(a,b,c,d) + \tfrac{1}{3} \widetilde{\aki}_6(a,c,b,d)
  \\
  &\quad
  -\tfrac{1}{3} \langle a,b,c,d \rangle - \tfrac{1}{3} \langle a,c,b,d \rangle + \tfrac{1}{6} \Phi_{1,3}(a,b,c,d),
  \end{align*}
where we have rewritten the Akivis elements in terms of the Sabinin operations:
  \begin{align*}
  \widetilde{\aki}_1(a,b,c,d) &= - \langle\langle\langle a, b \rangle, c \rangle, d \rangle,
  \\
  \widetilde{\aki}_2(a,b,c,d) &= - \langle \Phi_{1,2}( a, b, c ), d \rangle + \tfrac{1}{2} \langle \langle a, b, c \rangle, d \rangle ,
  \\
  \widetilde{\aki}_3(a,b,c,d) &= - \langle\langle a, b \rangle, \langle c, d \rangle\rangle,
  \\
  \widetilde{\aki}_4(a,b,c,d) &= - \Phi_{1,2}( \langle a, b \rangle, c , d ) + \tfrac{1}{2} \langle \langle a, b \rangle, c , d \rangle ,
  \\
  \widetilde{\aki}_5(a,b,c,d) &= - \Phi_{1,2}( a, \langle b , c \rangle , d ) + \tfrac{1}{2} \langle a , \langle b, c \rangle , d \rangle ,
  \\
  \widetilde{\aki}_6(a,b,c,d) &= - \Phi_{1,2}( a, b, \langle c , d \rangle ) + \tfrac{1}{2} \langle a, b , \langle c , d \rangle \rangle.
\end{align*}
Using the identities defining degree 4 Sabinin algebras, straightforward computations show that these operations
satisfy all the identities of a BTQQ algebra. For example, to prove the identity
  \begin{align*}
  [ b, (a,c,d) ] - [ b, (a,d,c) ] - ( a, b, [c,d] )
  &\equiv
  \lq a, b,c,d \rq - \lq a,b,d,c \rq
  \\
  &\quad
  - \lqq a,b,c,d \rqq + \lqq a,b,d,c \rqq,
  \end{align*}
we calculate as follows:
  \allowdisplaybreaks
  \begin{align*}
  &
  \lq a, b,c,d \rq
  - \lq a,b,d,c \rq
  - \lqq a,b,c,d \rqq
  + \lqq a,b,d,c \rqq
  \\
  &=
  - \widetilde{\aki}_2(a,c,d,b) + \tfrac{1}{4} \widetilde{\aki}_4(a,b,c,d)
  + \tfrac{1}{4} \widetilde{\aki}_4(a,b,d,c) - \tfrac{1}{2} \langle a,b,c,d \rangle
  \\
  &\quad
  + \tfrac{1}{4} \Phi_{2,2}(a,b,c,d) + \widetilde{\aki}_2(a,d,c,b)
  - \tfrac{1}{4} \widetilde{\aki}_4(a,b,d,c) - \tfrac{1}{4} \widetilde{\aki}_4(a,b,c,d)
  \\
  &\quad
  + \tfrac{1}{2} \langle a,b,d,c \rangle - \tfrac{1}{4} \Phi_{2,2}(a,b,d,c)
  - \tfrac{1}{3} \widetilde{\aki}_2(a,b,c,d) - \tfrac{2}{3} \widetilde{\aki}_2(a,b,d,c)
  \\
  &\quad
  - \tfrac{1}{3} \widetilde{\aki}_2(a,c,b,d) + \tfrac{1}{3} \widetilde{\aki}_2(a,c,d,b)
  + \tfrac{1}{2} \widetilde{\aki}_5(a,b,c,d) + \tfrac{1}{6} \widetilde{\aki}_5(a,b,d,c)
  \\
  &\quad
  + \tfrac{1}{6} \widetilde{\aki}_5(a,c,d,b) + \tfrac{1}{3} \widetilde{\aki}_6(a,b,c,d)
  + \tfrac{1}{3} \widetilde{\aki}_6(a,c,b,d) - \tfrac{1}{3} \langle a,b,c,d \rangle
  \\
  &\quad
  - \tfrac{1}{3} \langle a,c,b,d \rangle + \tfrac{1}{6} \Phi_{1,3}(a,b,c,d)
  - \tfrac{1}{3} \widetilde{\aki}_2(a,b,d,c) - \tfrac{2}{3} \widetilde{\aki}_2(a,b,c,d)
  \\
  &\quad
  - \tfrac{1}{3} \widetilde{\aki}_2(a,d,b,c) + \tfrac{1}{3} \widetilde{\aki}_2(a,d,c,b)
  + \tfrac{1}{2} \widetilde{\aki}_5(a,b,d,c) + \tfrac{1}{6} \widetilde{\aki}_5(a,b,c,d)
  \\
  &\quad
  + \tfrac{1}{6} \widetilde{\aki}_5(a,d,c,b) + \tfrac{1}{3} \widetilde{\aki}_6(a,b,d,c)
  + \tfrac{1}{3} \widetilde{\aki}_6(a,d,b,c) - \tfrac{1}{3} \langle a,b,d,c \rangle
  \\
  &\quad
  - \tfrac{1}{3} \langle a,d,b,c \rangle + \tfrac{1}{6} \Phi_{1,3}(a,b,d,c)
  \\
  &=
  - \widetilde{\aki}_2(a,c,d,b) + \widetilde{\aki}_2(a,d,c,b)
  - \langle a,b,c,d \rangle - \tfrac{1}{3} \widetilde{\aki}_2(a,b,c,d)
  \\
  &\quad
  + \tfrac{1}{3} \widetilde{\aki}_2(a,b,d,c) + \tfrac{1}{3} \widetilde{\aki}_2(a,c,b,d)
  - \tfrac{1}{3} \widetilde{\aki}_2(a,c,d,b) - \tfrac{1}{3} \widetilde{\aki}_2(a,d,b,c)
  \\
  &\quad
  + \tfrac{1}{3} \widetilde{\aki}_2(a,d,c,b) - \tfrac{1}{3} \widetilde{\aki}_5(a,b,c,d)
  + \tfrac{1}{3} \widetilde{\aki}_5(a,b,d,c) - \tfrac{1}{3} \widetilde{\aki}_5(a,c,d,b)
  \\
  &\quad
  - \tfrac{2}{3} \widetilde{\aki}_6(a,b,c,d) - \tfrac{1}{3} \widetilde{\aki}_6(a,c,b,d)
  + \tfrac{1}{3} \widetilde{\aki}_6(a,d,b,c) + \tfrac{2}{3} \langle a,b,c,d \rangle
  \\
  &\quad
  - \tfrac{1}{3} \langle a,c,b,d \rangle -\tfrac{1}{3} \langle a,d,b,c \rangle
  \\
  &=
  - \widetilde{\aki}_2(a,c,d,b) + \widetilde{\aki}_2(a,d,c,b)
  - \tfrac{1}{3} \big( \langle a,b,c,d \rangle + \langle a,c,d,b \rangle + \langle a,d,b,c \rangle \big)
  \\
  &\quad
  + \tfrac{1}{3} \big(
  \langle \Phi_{1,2}(a,b,c) , d \rangle - \tfrac{1}{2} \langle \langle a,b,c \rangle, d \rangle
  - \langle \Phi_{1,2}(a,b,d) , c \rangle + \tfrac{1}{2} \langle \langle a,b,d \rangle, c \rangle
  \\
  &\qquad
  - \langle \Phi_{1,2}(a,c,b) , d \rangle + \tfrac{1}{2} \langle \langle a,c,b \rangle, d \rangle
  + \langle \Phi_{1,2}(a,c,d) , b \rangle - \tfrac{1}{2} \langle \langle a,c,d \rangle, b \rangle
  \\
  &\qquad
  + \langle \Phi_{1,2}(a,d,b) , c \rangle - \tfrac{1}{2} \langle \langle a,d,b \rangle, c \rangle
  - \langle \Phi_{1,2}(a,d,c) , b \rangle + \tfrac{1}{2} \langle \langle a,d,c \rangle, b \rangle
  \big)
  \\
  &\quad
  + \tfrac{1}{3} \big(
  \Phi_{1,2}(a,\langle b,c \rangle,d) - \tfrac{1}{2} \langle a, \langle b,c \rangle , d \rangle
  - \Phi_{1,2}(a,\langle b,d \rangle,c) + \tfrac{1}{2} \langle a, \langle b,d \rangle , c \rangle
  \\
  &\qquad
  + \Phi_{1,2}(a,\langle c,d \rangle,b) - \tfrac{1}{2} \langle a, \langle c,d \rangle , b \rangle
  \big)
  \\
  &\quad
  + \tfrac{1}{3} \big(
  2 \Phi_{1,2}(a,b, \langle c,d \rangle) - \langle a,b \langle c,d \rangle \rangle
  + \Phi_{1,2}(a,c, \langle b,d \rangle) - \tfrac{1}{2} \langle a,c, \langle b,d \rangle \rangle
  \\
  &\qquad
  - \Phi_{1,2}(a,d, \langle b,c \rangle) + \tfrac{1}{2} \langle a,d, \langle b,c \rangle \rangle
  \big)
  \\
  &=
  - \widetilde{\aki}_2(a,c,d,b) + \widetilde{\aki}_2(a,d,c,b)
  \\
  &\quad
  + \tfrac{1}{3} \big(
  \langle \langle a,d,b \rangle , c \rangle + \langle a, \langle d,b \rangle , c \rangle + \langle \langle a,b,c \rangle , d \rangle
  \\
  &\qquad
  + \langle a, \langle b,c \rangle , d \rangle + \langle \langle a,c,d \rangle , b \rangle + \langle a, \langle c,d \rangle , b \rangle
  \big)
  \\
  &\quad
  + \tfrac{1}{3} \big(
  - \langle \langle a,b,c \rangle , d \rangle - \langle \langle a,c,d \rangle , b \rangle
  - \langle \langle a,d,b \rangle , c \rangle
  \big)
  \\
  &\quad
  + \tfrac{1}{3} \big(
  3\Phi_{1,2}(a,b,\langle c,d \rangle) - \langle a,\langle b,c \rangle,d \rangle
  -\langle a,\langle d,b \rangle, c \rangle - \tfrac{1}{2}\langle a,b,\langle c,d \rangle \rangle
  \big)
  \\
  &=
  -\widetilde{\aki}_2(a,c,d,b) + \widetilde{\aki}_2(a,d,c,b) - \widetilde{\aki}_6(a,b,c,d)
  \\
  &=
  -\aki_2(a,c,d,b) + \aki_2(a,d,c,b) - \aki_6(a,b,c,d).
  \end{align*}
The second equality is obtained by grouping terms, using the symmetries of $\Phi_{1,3}$ and $\Phi_{2,2}$
and the skew-symmetries of $\langle -, -, -, - \rangle$, $\widetilde{\aki}_5$ and $\widetilde{\aki}_6$.
For the third equality, we expand the Akivis elements $\widetilde{\aki}_i$.
The fourth equality follows from identity \eqref{S3} for degree 4 Sabinin algebras ($r=1$),
the symmetry of $\Phi_{1,2}$ and the skew-symmetry of $\langle -, - \rangle$ and $\langle -, - , - \rangle$.
The last two equalities are immediate from the definitions.
\end{proof}

%%%%%%%%%%%%%%%%%%%%%%%%%%%%%%%%%%%%%%%%%%%%%%%%%%%%%%%%%%%%%%%%%%%%%%%%%%%%%%%%%%%%

\section{Tangent algebras of monoassociative loops (BTQ algebras)}

According to Table \ref{algebratable}, the tangent algebras of monoassociative loops
may be defined in terms of the polynomial identities satisfied by the commutator, associator,
and first quaternator in a free power associative algebra.
By a result of Albert \cite{Albert}, we know that a nonassociative algebra over a field of characteristic 0
is power associative if and only if it satisfies the identities
  \[
  ( a, a, a ) \equiv 0,
  \qquad
  ( a^2, a, a ) \equiv 0,
  \]
which we call third-power and fourth-power associativity.

\begin{lemma} \label{degree3lemma}
In any third-power associative algebra, every multilinear polynomial identity of degree 3
satisfied by the commutator and associator is a consequence of
  \begin{align*}
  &
  (a,b,c) + (a,c,b) + (b,a,c) + (b,c,a) + (c,a,b) + (c,b,a) \equiv 0,
  \\
  &
  [[a,b],c] + [[b,c],a] + [[c,a],b] \equiv 2(a,b,c) + 2(b,c,a) + 2(c,a,b).
  \end{align*}
\end{lemma}

\begin{proof}
The first identity is the linearization of third-power associativity.
The second identity is the sum of the first identity and the Akivis identity.
To prove that every identity is a consequence of these two identities, we calculate as follows.
We construct a $10 \times 21$ matrix $E$ in which
the first 12 columns correspond to the multilinear nonassociative monomials,
  \[
  (ab)c, \; (ac)b, \; (ba)c, \; (bc)a, \; (ca)b, \; (cb)a, \;
  a(bc), \; a(cb), \; b(ac), \; b(ca), \; c(ab), \; c(ba),
  \]
and the last 9 columns correspond to the multilinear binary-ternary monomials,
  \[
  [[a,b],c], \; [[a,c],b], \; [[b,c],a], \;
  (a,b,c), \; (a,c,b), \; (b,a,c), \; (b,c,a), \; (c,a,b), \; (c,b,a).
  \]
The matrix has the following block structure,
  \[
  E = \left[ \begin{array}{c|c} T & O \\ \midrule X & I \end{array} \right],
  \]
where $T$ is the $1 \times 12$ vector of coefficients of the third-power associativity,
$O$ is the $1 \times 9$ zero vector, $X$ is the $9 \times 12$ matrix in which the $(i,j)$ entry
is the coefficient of the $j$-th nonassociative monomial in the expansion of the $i$-th
binary-ternary monomial, and $I$ is the $9 \times 9$ identity matrix. (See Table \ref{degree3}.)
\break
\begin{table}[h]
\[
\left[
\begin{array}{rrrrrrrrrrrr|rrrrrrrrr}
1 &\!\! 1 &\!\! 1 &\!\! 1 &\!\! 1 &\!\! 1 &\!\! -1 &\!\! -1 &\!\! -1 &\!\! -1 &\!\! -1 &\!\! -1 & 0 &\!\! 0 &\!\! 0 &\!\! 0 &\!\! 0 &\!\! 0 &\!\! 0 &\!\! 0 &\!\! 0 \\ \midrule
1 &\!\! 0 &\!\! -1 &\!\! 0 &\!\! 0 &\!\! 0 &\!\! 0 &\!\! 0 &\!\! 0 &\!\! 0 &\!\! -1 &\!\! 1 & 1 &\!\! 0 &\!\! 0 &\!\! 0 &\!\! 0 &\!\! 0 &\!\! 0 &\!\! 0 &\!\! 0 \\
0 &\!\! 1 &\!\! 0 &\!\! 0 &\!\! -1 &\!\! 0 &\!\! 0 &\!\! 0 &\!\! -1 &\!\! 1 &\!\! 0 &\!\! 0 & 0 &\!\! 1 &\!\! 0 &\!\! 0 &\!\! 0 &\!\! 0 &\!\! 0 &\!\! 0 &\!\! 0 \\
0 &\!\! 0 &\!\! 0 &\!\! 1 &\!\! 0 &\!\! -1 &\!\! -1 &\!\! 1 &\!\! 0 &\!\! 0 &\!\! 0 &\!\! 0 & 0 &\!\! 0 &\!\! 1 &\!\! 0 &\!\! 0 &\!\! 0 &\!\! 0 &\!\! 0 &\!\! 0 \\
1 &\!\! 0 &\!\! 0 &\!\! 0 &\!\! 0 &\!\! 0 &\!\! -1 &\!\! 0 &\!\! 0 &\!\! 0 &\!\! 0 &\!\! 0 & 0 &\!\! 0 &\!\! 0 &\!\! 1 &\!\! 0 &\!\! 0 &\!\! 0 &\!\! 0 &\!\! 0 \\
0 &\!\! 1 &\!\! 0 &\!\! 0 &\!\! 0 &\!\! 0 &\!\! 0 &\!\! -1 &\!\! 0 &\!\! 0 &\!\! 0 &\!\! 0 & 0 &\!\! 0 &\!\! 0 &\!\! 0 &\!\! 1 &\!\! 0 &\!\! 0 &\!\! 0 &\!\! 0 \\
0 &\!\! 0 &\!\! 1 &\!\! 0 &\!\! 0 &\!\! 0 &\!\! 0 &\!\! 0 &\!\! -1 &\!\! 0 &\!\! 0 &\!\! 0 & 0 &\!\! 0 &\!\! 0 &\!\! 0 &\!\! 0 &\!\! 1 &\!\! 0 &\!\! 0 &\!\! 0 \\
0 &\!\! 0 &\!\! 0 &\!\! 1 &\!\! 0 &\!\! 0 &\!\! 0 &\!\! 0 &\!\! 0 &\!\! -1 &\!\! 0 &\!\! 0 & 0 &\!\! 0 &\!\! 0 &\!\! 0 &\!\! 0 &\!\! 0 &\!\! 1 &\!\! 0 &\!\! 0 \\
0 &\!\! 0 &\!\! 0 &\!\! 0 &\!\! 1 &\!\! 0 &\!\! 0 &\!\! 0 &\!\! 0 &\!\! 0 &\!\! -1 &\!\! 0 & 0 &\!\! 0 &\!\! 0 &\!\! 0 &\!\! 0 &\!\! 0 &\!\! 0 &\!\! 1 &\!\! 0 \\
0 &\!\! 0 &\!\! 0 &\!\! 0 &\!\! 0 &\!\! 1 &\!\! 0 &\!\! 0 &\!\! 0 &\!\! 0 &\!\! 0 &\!\! -1 & 0 &\!\! 0 &\!\! 0 &\!\! 0 &\!\! 0 &\!\! 0 &\!\! 0 &\!\! 0 &\!\! 1
\end{array}
\right]
\]
\caption{Matrix for proof of Lemma \ref{degree3lemma}}
\label{degree3}
\end{table}

We compute the row canonical form (RCF) of this matrix, and identify the rows whose leading 1s are in the right part
of the matrix. (See Table \ref{degree3rcf}.)
\break
\begin{table}[h]
\[
\left[
\begin{array}{rrrrrrrrrrrr|rrrrrrrrr}
1 &\!\! 0 &\!\! 0 &\!\! 0 &\!\! 0 &\!\! 0 &\!\! 0 &\!\! 0 &\!\! -1 &\!\! 0 &\!\! -1 &\!\! 1 & 0 &\!\! 1 &\!\! -1 &\!\! 0 &\!\! -2 &\!\! -1 &\!\! 0 &\!\! 0 &\!\! -2
\\
0 &\!\! 1 &\!\! 0 &\!\! 0 &\!\! 0 &\!\! 0 &\!\! 0 &\!\! 0 &\!\! -1 &\!\! 1 &\!\! -1 &\!\! 0 & 0 &\!\! 1 &\!\! 0 &\!\! 0 &\!\! 0 &\!\! 0 &\!\! 0 &\!\! 1 &\!\! 0
\\
0 &\!\! 0 &\!\! 1 &\!\! 0 &\!\! 0 &\!\! 0 &\!\! 0 &\!\! 0 &\!\! -1 &\!\! 0 &\!\! 0 &\!\! 0 & 0 &\!\! 0 &\!\! 0 &\!\! 0 &\!\! 0 &\!\! 1 &\!\! 0 &\!\! 0 &\!\! 0
\\
0 &\!\! 0 &\!\! 0 &\!\! 1 &\!\! 0 &\!\! 0 &\!\! 0 &\!\! 0 &\!\! 0 &\!\! -1 &\!\! 0 &\!\! 0 & 0 &\!\! 0 &\!\! 0 &\!\! 0 &\!\! 0 &\!\! 0 &\!\! 1 &\!\! 0 &\!\! 0
\\
0 &\!\! 0 &\!\! 0 &\!\! 0 &\!\! 1 &\!\! 0 &\!\! 0 &\!\! 0 &\!\! 0 &\!\! 0 &\!\! -1 &\!\! 0 & 0 &\!\! 0 &\!\! 0 &\!\! 0 &\!\! 0 &\!\! 0 &\!\! 0 &\!\! 1 &\!\! 0
\\
0 &\!\! 0 &\!\! 0 &\!\! 0 &\!\! 0 &\!\! 1 &\!\! 0 &\!\! 0 &\!\! 0 &\!\! 0 &\!\! 0 &\!\! -1 & 0 &\!\! 0 &\!\! 0 &\!\! 0 &\!\! 0 &\!\! 0 &\!\! 0 &\!\! 0 &\!\! 1
\\
0 &\!\! 0 &\!\! 0 &\!\! 0 &\!\! 0 &\!\! 0 &\!\! 1 &\!\! 0 &\!\! -1 &\!\! 0 &\!\! -1 &\!\! 1 & 0 &\!\! 1 &\!\! -1 &\!\! 0 &\!\! -1 &\!\! 0 &\!\! 1 &\!\! 1 &\!\! -1
\\
0 &\!\! 0 &\!\! 0 &\!\! 0 &\!\! 0 &\!\! 0 &\!\! 0 &\!\! 1 &\!\! -1 &\!\! 1 &\!\! -1 &\!\! 0 & 0 &\!\! 1 &\!\! 0 &\!\! 0 &\!\! -1 &\!\! 0 &\!\! 0 &\!\! 1 &\!\! 0
\\ \midrule
0 &\!\! 0 &\!\! 0 &\!\! 0 &\!\! 0 &\!\! 0 &\!\! 0 &\!\! 0 &\!\! 0 &\!\! 0 &\!\! 0 &\!\! 0 & 1 &\!\! -1 &\!\! 1 &\!\! 0 &\!\! 2 &\!\! 2 &\!\! 0 &\!\! 0 &\!\! 2
\\
0 &\!\! 0 &\!\! 0 &\!\! 0 &\!\! 0 &\!\! 0 &\!\! 0 &\!\! 0 &\!\! 0 &\!\! 0 &\!\! 0 &\!\! 0 & 0 &\!\! 0 &\!\! 0 &\!\! 1 &\!\! 1 &\!\! 1 &\!\! 1 &\!\! 1 &\!\! 1
\end{array}
\right]
\]
\caption{RCF of matrix from Table \ref{degree3}}
\label{degree3rcf}
\end{table}

The last two rows of the RCF contain the coefficient vectors of the two identities in the statement of this Lemma.
\end{proof}

\begin{definition} \label{deflifting}
From a multilinear polynomial identity $f( a_1, a_2, \dots, a_n )$ in degree $n$ for a nonassociative algebra 
with one binary operation, we obtain $n+2$ identities in degree $n+1$ which generate all the consequences
of $f$ in degree $n+1$ as a module over the symmetric group.  We introduce a new variable $a_{n+1}$ and perform
$n$ substitutions and two multiplications, obtaining the \textbf{liftings} of $f$:
  \begin{align*}
  &
  f( a_1 a_{n+1}, a_2, \dots, a_n ),
  \quad
  f( a_1, a_2 a_{n+1}, \dots, a_n ),
  \quad
  \dots,
  \quad
  f( a_1, a_2, \dots, a_n a_{n+1} ),
  \\
  &
  f( a_1, a_2, \dots, a_n ) a_{n+1},
  \quad
  a_{n+1} f( a_1, a_2, \dots, a_n ).
  \end{align*}
\end{definition}

\begin{lemma} \label{tpalemma}
Every consequence in degree 4 of third-power associativity
is a linear combination of permutations of these three identities:
  \[
  \begin{array}{l}
  \tpa_1(a,b,c,d) =
  ( ad, b, c ) + ( ad, c, b ) + ( b, ad, c )
  \\
  \qquad\qquad\qquad\qquad
  +
  ( b, c, ad ) + ( c, ad, b ) + ( c, b, ad ),
  \\[4pt]
  \tpa_2(a,b,c,d) =
  ( a, b, c ) d + ( a, c, b ) d + ( b, a, c ) d
  \\
  \qquad\qquad\qquad\qquad
  +
  ( b, c, a ) d + ( c, a, b ) d + ( c, b, a ) d,
  \\[4pt]
  \tpa_3(a,b,c,d) =
  d ( a, b, c ) + d ( a, c, b ) + d ( b, a, c )
  \\
  \qquad\qquad\qquad\qquad
  +
  d ( b, c, a ) + d ( c, a, b ) + d ( c, b, a ).
  \end{array}
  \]
\end{lemma}

\begin{proof}
When lifting the third-power associativity to degree 4, symmetry implies that
it suffices to consider one substitution and two multiplications\footnote{We note that
Bremner et al.~\cite{BHPTU} contains an error: the two equations at the top of page 173
should be three equations, and the commutators $[a,b]$ should be nonassociative products $ab$.}.
\end{proof}

\begin{lemma}
In any third-power associative algebra, the second quaternator is redundant:
$\qua_2$ is a linear combination of permutations of $\qua_1$,
the Akivis elements, and the consequences in degree 4 of third-power associativity.
\end{lemma}

\begin{proof}
We construct a matrix with a $120 \times 120$ upper block and a $24 \times 120$ lower block.
The columns correspond to the multilinear nonassociative monomials in degree 4: each of the five association types
allows 24 permutations of $a,b,c,d$.
We process the following ten polynomials in order: the three identities of Lemma \ref{tpalemma}, the six Akivis elements \eqref{akivis},
and the two quaternators \eqref{quaternators}.  At each step we put the coefficient vectors of the permutations
of the polynomial in the lower block, and then compute the row canonical form.  We obtain the following sequence of
ranks: 12, 16, 20, 32, 46, 49, 61, 73, 73, 82, 82.  In particular, the second quaternator belongs to the $S_4$-submodule
generated by the other nine elements.
\end{proof}

\begin{remark}
To obtain an explicit formula, we construct a matrix consisting of ten $120 \times 24$ blocks from left to right.
The rows correspond to the multilinear nonassociative monomials in degree 4.  In the columns of each block, we store
the coefficient vectors of all permutations of the first quaternator, the six Akivis elements, and the three identities
of Lemma \ref{tpalemma}.  We compute the RCF and find that the rank is 82.  We identify the columns of the RCF which
contain the leading 1s of the nonzero rows; the corresponding columns of the original matrix form a basis of the column space.
We then construct a $120 \times 83$ matrix; in columns 1--82 we put the basis vectors of the column space
from the previous matrix, and in column 83 we put the coefficient vector of the second quaternator \eqref{quaternators}.
We compute the RCF; the last column contains the coefficients of the second
quaternator with respect to the basis of the column space.  This expression for $\qua_2(a,b,c,d)$ has 54 terms:
  \begin{align*}
  \tfrac14
  &\big(
   3 \qua_1 (a,b,c,d)
  -  \qua_1 (a,b,d,c)
  +  \qua_1 (a,c,b,d)
  -  \qua_1 (a,c,d,b)
  +3 \qua_1 (a,d,b,c)
  \\
  &
  -  \qua_1 (a,d,c,b)
  -  \qua_1 (b,c,a,d)
  -3 \qua_1 (b,d,a,c)
  +2 \qua_1 (c,a,d,b)
  +  \qua_1 (c,b,a,d)
  \\
  &
  -2 \qua_1 (c,b,d,a)
  -  \qua_1 (c,d,a,b)
  -2 \qua_1 (d,a,b,c)
  +2 \qua_1 (d,a,c,b)
  +3 \qua_1 (d,b,a,c)
  \\
  &
  -2 \qua_1 (d,b,c,a)
  -  \qua_1 (d,c,a,b)
  -2 \aki_1 (a,d,b,c)
  +2 \aki_1 (a,d,c,b)
  -  \aki_1 (b,c,d,a)
  \\
  &
  +2 \aki_1 (b,d,a,c)
  -  \aki_1 (b,d,c,a)
  -  \aki_1 (c,d,b,a)
  +2 \aki_2 (a,b,c,d)
  -2 \aki_2 (a,c,b,d)
  \\
  &
  +4 \aki_2 (a,c,d,b)
  -2 \aki_2 (b,a,c,d)
  +2 \aki_2 (b,c,a,d)
  -2 \aki_2 (b,d,c,a)
  +2 \aki_2 (c,d,b,a)
  \\
  &
  +4 \aki_2 (d,a,c,b)
  -2 \aki_2 (d,b,c,a)
  +2 \aki_3 (a,c,b,d)
  +2 \aki_3 (a,d,b,c)
  +2 \aki_5 (a,b,c,d)
  \\
  &
  +2 \aki_5 (a,b,d,c)
  -  \aki_5 (b,a,c,d)
  -  \aki_5 (b,a,d,c)
  +  \aki_5 (c,a,b,d)
  +  \aki_5 (c,a,d,b)
  \\
  &
  -2 \aki_5 (c,b,d,a)
  +  \aki_5 (d,a,b,c)
  -  \aki_5 (d,a,c,b)
  +2 \aki_6 (a,b,c,d)
  -2 \aki_6 (a,c,b,d)
  \\
  &
  +2 \aki_6 (b,c,a,d)
  -  \tpa_1 (a,b,c,d)
  +  \tpa_1 (a,b,d,c)
  -  \tpa_1 (a,c,d,b)
  +2 \tpa_1 (b,a,c,d)
  \\
  &
  -  \tpa_2 (a,b,c,d)
  -  \tpa_2 (a,b,d,c)
  -  \tpa_2 (a,c,d,b)
  +  \tpa_2 (b,c,d,a)
  \big).
  \end{align*}
\end{remark}

The next theorem gives the polynomial identities of degree 4 relating commutator, associator,
and first quaternator in any power associative algebra.  These identities define a new variety
of nonassociative multioperator algebras generalizing Lie, Malcev, and Bol algebras.

\begin{theorem} \label{degree4theorem}
In any power associative algebra, every multilinear polynomial identity of degree 4
satisfied by the commutator, associator, and first quaternator, is a consequence of
the skew-symmetry of the commutator, the two identities of Lemma \ref{degree3lemma},
and these four identities:
  \allowdisplaybreaks
  \begin{align*}
  &
    [(a,b,c),d]
  - [(d,b,c),a]
  - ([a,d],b,c)
  + \lq a,d,b,c \rq
  - \lq d,a,b,c \rq
  \equiv 0,
  \\[4pt]
  &
    (a,[b,c],d)
  - (a,[b,d],c)
  + (a,[c,d],b)
  - (a,b,[c,d])
  + (a,c,[b,d])
  - (a,d,[b,c])
  \\
  &
  - \lq a,b,c,d \rq
  + \lq a,b,d,c \rq
  + \lq a,c,b,d \rq
  - \lq a,c,d,b \rq
  - \lq a,d,b,c \rq
  + \lq a,d,c,b \rq
  \equiv 0,
  \\[4pt]
  &
    [(a,b,c),d]
  - [(a,d,c),b]
  - [(c,b,a),d]
  + [(c,d,a),b]
  + [(b,a,d),c]
  - [(b,c,d),a]
  \\
  &
  - [(d,a,b),c]
  + [(d,c,b),a]
  - ([a,b],c,d)
  - ([a,d],b,c)
  - ([c,b],d,a)
  - ([c,d],a,b)
  \\
  &
  - (a,[c,b],d)
  - (c,[a,d],b)
  - (b,[c,d],a)
  - (d,[a,b],c)
  - (a,b,[c,d])
  - (c,d,[a,b])
  \\
  &
  - (b,c,[a,d])
  - (d,a,[c,b])
  \equiv 0,
  \\[4pt]
  &
  \sum_{\pi} \lq a^\pi, b^\pi, c^\pi, d^\pi \rq
  \equiv 0,
  \end{align*}
where $\pi$ runs over all permutations of $\{a,b,c,d\}$.
\end{theorem}

\begin{proof}
We first show that these identities are satisfied by the three operations in every
power associative algebra.

The first identity is equivalent to the third identity of Theorem \ref{BTQQtheorem}.
To verify the second identity, we expand it in the free nonassociative algebra
using the definitions of the commutator, associator, and first quaternator,
and observe that the result collapses to zero.  In other words, this identity
is a consequence of the first three identities of Theorem \ref{BTQQtheorem}.
The third identity can be verified in the same way, but it is an Akivis element
(it involves only the commutator and associator), so it is a consequence of the
first two identities of Theorem \ref{BTQQtheorem}.

To prove the fourth identity, we expand the quaternators but not the commutators
and associators.  We obtain the following expression:
  \[
  \fpa(a,b,c,d)
  - \sigma_{a,b,c,d} \tpa_2(a,b,c,d)
  - \sigma_{a,b,c,d} \tpa_3(a,b,c,d),
  \]
where $\fpa(a,b,c,d)$ is the linearization of the fourth-power associativity,
and $\sigma_{a,b,c,d}$ denotes the cyclic sum over the four arguments.  This expression
obviously vanishes in every power associative algebra.

To complete the proof, we need to show that every identity satisfied by the three
operations in every power associative algebra is a consequence of the identities
stated in the theorem.  For this we use computer algebra, following the basic ideas
in the proof of Lemma \ref{degree3lemma}.

There are seven possible combinations of the three operations in degree 4:
  \begin{alignat*}{4}
  &[ [ [ a, b ], c ], d ], &\qquad
  &[ ( a, b, c ), d ], &\qquad
  &[ [ a, b ], [ c, d ] ], &\qquad
  &( [ a, b ], c, d ),
  \\
  &( a, [ b, c ], d ), &\qquad
  &( a, b, [ c, d ] ), &\qquad
  &\lq a, b, c, d \rq.
  \end{alignat*}
Owing to the skew-symmetry of the commutator, these elements admit respectively
12, 24, 3, 12, 12, 12, 24 distinct permutations, for a total of 99 multilinear monomials,
which we call BTQ monomials (since we now have only the first quaternator).

A power associative algebra satisfies these multilinear identities in degree 4:
the liftings of third-power associativity, $\tpa_i(a,b,c,d)$ for $i = 1, 2, 3$,
and the linearized fourth-power associativity $\fpa(a,b,c,d)$.
By symmetries of these identities, a linear basis for the submodule generated by
these identities consists of all 24 permutations of $\tpa_1$,
the four cyclic permutations of $\tpa_2$ and $\tpa_3$, and
the identity permutation of $\fpa$.

We now construct a block matrix as in the proof of Lemma \ref{degree3lemma}
but larger:
  \[
  E = \left[ \begin{array}{c|c} T\!F & O \\ \midrule X & I \end{array} \right].
  \]
The left part of the matrix has 120 columns corresponding to the multilinear nonassociative monomials.
The right part has 99 columns corresponding to the multilinear BTQ monomials.
Block $T\!F$ contains the multilinear identities in degree 4 for a power associative algebra.
Block $X$ contains the expansions of the multilinear BTQ monomials into the free nonassociative algebra
using the definitions of commutator, associator, and first quaternator.
Blocks $I$ and $O$ are the identity matrix and the zero matrix.

All the entries of $E$ are integers, so we may compute the Hermite normal form (HNF) of $E$ over $\mathbb{Z}$.
The rank is 120, and the first 82 rows of the HNF have leading 1s in columns 1--120.
The remaining 38 rows have leading 1s in columns 121--219.
These rows form a lattice basis for the polynomial identities with integer coefficients
satisfied by the three operations in every power associative algebra.
We now apply the LLL algorithm for lattice basis reduction to these 38 vectors, and
sort the reduced vectors by increasing Euclidean length.
(For the application of lattice algorithms to polynomial identities, see \cite{BP}.
For a general introduction to lattice basis reduction, see \cite{B}.)

We generate all liftings of the two identities from Lemma \ref{degree3lemma} to degree 4;
in this case the binary operation is the commutator, and so the last lifting in
Definition \ref{deflifting} is redundant.
Following the same basic idea as in the proof of Lemma \ref{degree3lemma}, we determine a basis
for the space of identities in degree 4 which follow from the identities in degree 3.  We then
process the identities resulting from the reduced lattice basis to determine a set of generators
for the new identities in degree 4.
These generators are the identities stated in the theorem.
\end{proof}

\begin{definition} \label{btqdefinition}
A \textbf{BTQ algebra} is a vector space with an anticommutative bilinear operation $[a,b]$,
a trilinear operation $(a,b,c)$, and a quadrilinear operation $\lq a,b,c,d \rq$, satisfying
the identities of Lemma \ref{degree3lemma} and Theorem \ref{degree4theorem}.
\end{definition}

%%%%%%%%%%%%%%%%%%%%%%%%%%%%%%%%%%%%%%%%%%%%%%%%%%%%%%%%%%%%%%%%%%%%%%%%%%%%%%%%%%%%

\section{Special identities for BTQ algebras}

In this section we show that there exist identities of degrees 5 and 6 which are satisfied by
the commutator, associator, and first quaternator in every power associative algebra,
but which are not consequences of the defining identities for BTQ algebras from
Definition \ref{btqdefinition}.  We present an explicit identity in degree 5, but
only a computational existence proof for degree 6.  In order to reduce the amount of
computer memory that we need, we use the Wedderburn decomposition of the group algebra
$\mathbb{Q} S_n$ as direct sum of simple ideals which are isomorphic to full matrix
algebras; for general discussions of this method, see \cite{BP0,BP2,BP3}.

\begin{theorem} \label{newidentity}
The following nonlinear identity is satisfied by the commutator, associator, and
first quaternator in every power associative algebra, but is not a consequence of
the defining identities for BTQ algebras:
  \[
  [ ( a, [a, b], a ), a ] - [ \lq a, a, a, b \rq, a ] + [ \lq a, b, a, a \rq, a ]
  \equiv
  0.
  \]
\end{theorem}

\begin{proof}
There are seven distinct irreducible representations of the symmetric group $S_5$,
with dimensions $d = 1$, 4, 5, 6, 5, 4, 1 corresponding to partitions $\lambda = 5$,
41, 32, 311, 221, 2111, 11111.  The representation matrices $\rho_i(\pi)$ ($1 \le i \le 7$)
are given by the projections onto the corresponding simple ideals in the direct sum
decomposition of the group algebra:
  \[
  \mathbb{Q} S_5
  =
  \mathbb{Q} \oplus
  M_4(\mathbb{Q}) \oplus
  M_5(\mathbb{Q}) \oplus
  M_6(\mathbb{Q}) \oplus
  M_5(\mathbb{Q}) \oplus
  M_4(\mathbb{Q}) \oplus
  \mathbb{Q}.
  \]
These representation matrices may be computed following the method in \cite{Clifton}.
The basic idea of the following proof is to decompose a multilinear polynomial identity
into irreducible components corresponding to these simple matrix algebras.

Given any multilinear polynomial of degree 5 in the free nonassociative algebra,
we first sort its terms into the 14 possible association types,
  \begin{alignat*}{3}
  &( ( ( ( -, - ), - ), - ), - ), &\qquad
  &( ( ( -, ( -, - ) ), - ), - ), &\qquad
  &( ( ( -, - ), ( -, - ) ), - ),
  \\
  &( ( -, ( ( -, - ), - ) ), - ), &\qquad
  &( ( -, ( -, ( -, - ) ) ), - ), &\qquad
  &( ( ( -, - ), - ), ( -, - ) ),
  \\
  &( ( -, ( -, - ) ), ( -, - ) ), &\qquad
  &( ( -, - ), ( ( -, - ), - ) ), &\qquad
  &( ( -, - ), ( -, ( -, - ) ) ),
  \\
  &( -, ( ( ( -, - ), - ), - ) ), &\qquad
  &( -, ( ( -, ( -, - ) ), - ) ), &\qquad
  &( -, ( ( -, - ), ( -, - ) ) ),
  \\
  &( -, ( -, ( ( -, - ), - ) ) ), &\qquad
  &( -, ( -, ( -, ( -, - ) ) ) ).
  \end{alignat*}
Within each association type, the terms differ only by permutation of the variables,
and so the polynomial can be regarded as an element of
the direct sum of 14 copies of the group algebra $\mathbb{Q} S_5$.  For each
partition $\lambda$, corresponding to an irreducible representation of dimension $d$,
we compute the representation matrices of the 14 components of the polynomial, and
obtain a matrix of size $d \times 14d$.

We apply this to the consequences in degree 5 of the power associative identities
in the free nonassociative algebra.
Using Definition \ref{deflifting} we first lift the linearized third-power associativity to degree 4
using the nonassociative binary product, 
and then lift these identities and the linearized fourth-power
associativity to degree 5.  Altogether we obtain 36 multilinear identities
$P_i(a,b,c,d,e)$ ($1 \le i \le 36$) in degree 5 which are elements of the free
nonassociative algebra on 5 generators; all the multilinear identities in degree 5
satisfied by every power associative algebra are consequences of these identities.
Combining the representation matrices for these identities, we obtain a matrix $P$ of
size $36d \times 14d$.

We next consider the 22 distinct BTQ monomials with the identity permutation of the
variables:
  \begin{alignat*}{5}
  &[ [ [ [ a, b ], c ], d ], e ], &\;\;
  &[ [ ( a, b, c ), d ], e ], &\;\;
  &[ [ [ a, b ], [ c, d ] ], e ], &\;\;
  &[ ( [ a, b ], c, d ), e ], &\;\;
  &[ ( a, [ b, c ], d ), e ],
  \\
  &[ ( a, b, [ c, d ] ), e ], &\;\;
  &[ \lq a, b, c, d \rq, e ], &\;\;
  &[ [ [ a, b ], c ], [ d, e ] ], &\;\;
  &[ ( a, b, c ), [ d, e ] ], &\;\;
  &( [ [ a, b ], c ], d, e ),
  \\
  &( ( a, b, c ), d, e ), &\;\;
  &( [ a, b ], [ c, d ], e ), &\;\;
  &( [ a, b ], c, [ d, e ] ), &\;\;
  &( a, [ [ b, c ], d ], e ), &\;\;
  &( a, ( b, c, d ), e ),
  \\
  &( a, [ b, c ], [ d, e ] ), &\;\;
  &( a, b, [ [ c, d ], e ] ), &\;\;
  &( a, b, ( c, d, e ) ), &\;\;
  &\lq [ a, b ], c, d, e \rq, &\;\;
  &\lq a, [ b, c ], d, e \rq,
  \\
  &\lq a, b, [ c, d ], e \rq, &\;\;
  &\lq a, b, c, [ d, e ] \rq.
  \end{alignat*}
We write $X_i(a,b,c,d,e)$ ($1 \le i \le 22$) for the expansions of these monomials
in the free nonassociative algebra using the definitions of the commutator, associator
and first quaternator.
Combining the representation matrices for these expansions, we obtain a matrix $X$ of
size $22d \times 14d$.

We create a block matrix with the following structure,
  \[
  \left[
  \begin{array}{c|c}
  P & O \\
  \midrule
  X & I
  \end{array}
  \right],
  \]
where $O$ is the zero matrix and $I$ is the identity matrix.
We compute the row canonical form (RCF) of this matrix and identify the rows
whose leading 1s occur to the right of the vertical line.
These rows represent polynomial identities satisfied by the BTQ operations; more
precisely, they represent dependence relations among the BTQ monomials resulting
from dependence relations among their expansions in the free power associative algebra.
The number of these rows is displayed in column ``all identities'' of Table \ref{btq5ranks}.

We need to compare these identities with the skew-symmetries of the BTQ monomials
resulting from the anticommutativity of the commutator, and the liftings of the BTQ
identities from lower degrees.
Each BTQ monomial may have 0, 1 or more skew-symmetries depending on the occurrences
of the commutator.  For example, the monomial $([a,b],c,[d,e])$ has two skew-symmetries
that can be written as $I + T$ in the group algebra where $I$ is the identity permutation,
and $T$ is $bacde$ or $abced$.

The two BTQ identities in degree 3 (see Lemma \ref{degree3lemma}) can be lifted to degree 5
in two inequivalent ways: either using the commutator twice, or the associator once.
Lifting using the commutator has already been mentioned (see the last paragraph of the proof
of Theorem \ref{degree4theorem}); using a ternary operation the identity $I(a,b,c)$ in degree 3
produces six identities in degree 5:
  \begin{alignat*}{3}
  &I((a,d,e),b,c), &\qquad
  &I(a,(b,d,e),c), &\qquad
  &I(a,b,(c,d,e)), \\
  &(I(a,b,c),d,e), &\qquad
  &(d,I(a,b,c),e), &\qquad
  &(d,e,I(a,b,c)).
  \end{alignat*}
We also lift the four BTQ identities of degree 4 (see Theorem \ref{degree4theorem}) to degree 5
using the commutator.

For each partition $\lambda$ with corresponding representation of dimension $d$, we construct
a matrix with $22d$ columns containing the representation matrices for the symmetries and
liftings, and compute its RCF.
The rank of this matrix is displayed in column ``symm+lift'' of Table \ref{btq5ranks}.

\begin{table}[h]
\begin{tabular}{ccccc}
partition & dimension & symm+lift & all identities & new \\
\toprule
5 & 1 & 22 & 22 & 0 \\
41 & 4 & 73 & 74 & 1 \\
32 & 5 & 79 & 79 & 0 \\
311 & 6 & 93 & 93 & 0 \\
221 & 5 & 72 & 72 & 0 \\
2111 & 4 & 56 & 56 & 0 \\
11111 & 1 & 14 & 14 & 0 \\
\bottomrule
\end{tabular}
\medskip
\caption{Ranks for BTQ identities in degree 5}
\label{btq5ranks}
\end{table}

In Table \ref{btq5ranks}, column ``new'' gives the difference between columns
``all identities'' and ``symm+lift''.  When this value is zero, we also checked
that the two matrices in RCF resulting from the two computations are exactly
equal, and hence every BTQ identity for these partitions is a consequence of
the symmetries and liftings.  The only nonzero value occurs for partition 41,
and this indicates the existence of a new identity.

Since the new identity exists in partition 41, we expect to find a nonlinear
identity in two variables, of degree 4 in one variable $a$ and linear in the
other variable $b$.  There are 22 multilinear BTQ monomials with the identity
permutation, and hence 110 possible substitutions with some permutation of
$aaaab$.  Many of these are zero as a result of the identities $[a,a] \equiv 0$,
$(a,a,a) \equiv 0$ and $\lq a,a,a,a \rq \equiv 0$.  Eliminating these we are
left with 49 nonlinear BTQ monomials.  We now follow the same approach as in
the proof of Theorem \ref{degree4theorem} but using nonlinear monomials.
We compute the Hermite normal form of the matrix, and extract the lower right
block whose rows represent identities in the BTQ monomials.  We divide each
row by the GCD of its entries, and then apply lattice basis reduction to the
lattice spanned by the rows.  All except one of the elements of the reduced
lattice basis represent identities which are consequences of identities of
lower degree; the exception is the identity in the statement of this Theorem.

We checked this result by linearizing this new identity, combining the result
with the symmetries and lifted identities, and recomputing the representation
matrices.  The matrices we obtain from this calculation are now exactly the same
as the matrices for ``all identities'' in every partition.
\end{proof}

\begin{remark}
A direct proof that this new identity is satisfied in every power associative algebra
is given by the following formula, where
  \[
  T(x,y,z) = \sum_{\sigma} (x^\sigma, y^\sigma, z^\sigma ),
  \qquad
  F(w,x,y,z) = \sum_{\tau} ( w^\tau x^\tau, y^\tau, z^\tau ),
  \]
denote the multilinear forms of third- and fourth-power associativity.
We have
  \begin{align*}
  &
  [ ( a, [a, b], a ), a ] - [ \lq a, a, a, b \rq, a ] + [ \lq a, b, a, a \rq, a ]
  =
  \\
  \tfrac{1}{30}
  \big( \,
  &
    3  T((ab)a,a,a)
  -12  T((aa)b,a,a)
   +6  T((aa)a,a,b)
   +3  T((ba)a,a,a)
   \\
   &
   +6  T(aa,ab,a)
  -12  T(aa,aa,b)
   +6  T(aa,ba,a)
   +3  T(a(ab),a,a)
   \\
   &
   +3  T(a(ba),a,a)
   +6  T(a(aa),a,b)
  -12  T(b(aa),a,a)
   +3  T(aa,a,a)b
   \\
   &
  -12  T(ab,a,a)a
   +6  T(aa,a,b)a
   +3  T(ba,a,a)a
   -3  bT(aa,a,a)
   \\
   &
  +12  aT(ab,a,a)
   -6  aT(aa,a,b)
   -3  aT(ba,a,a)
   +2  (T(a,a,a)a)b
   \\
   &
   +7  (T(a,a,a)b)a
   -9  (T(a,a,b)a)a
   -2  b(T(a,a,a)a)
   -7  a(T(a,a,a)b)
   \\
   &
   +9  a(T(a,a,b)a)
   -   F(a,a,a,a)b
   +   F(a,a,a,b)a
   +   bF(a,a,a,a)
   \\
   &
   -   aF(a,a,a,b)
  \, \big).
  \end{align*}
\end{remark}

\begin{theorem}
There exists a nonlinear identity of degree 6 in two variables, with degree 3 in each variable,
which is satisfied by the commutator, associator, and first quaternator in every power associative algebra,
but is not a consequence of the defining identities for BTQ algebras and the new identity of Theorem \ref{newidentity}.
\end{theorem}

\begin{proof}
We extended the calculations in the proof of Theorem  \ref{newidentity} to degree 6 and obtained the results in Table \ref{btq6ranks}.
\end{proof}

\begin{table}[h]
\begin{tabular}{ccccc}
partition & dimension & symm+lift & all identities & new \\
\toprule
6 & 1 & 80 & 80 & 0 \\
51 & 5 & 367 & 367 & 0 \\
42 & 9 & 602 & 602 & 0 \\
411 & 10 & 662 & 662 & 0 \\
33 & 5 & 314 & 315 & 1 \\
321 & 16 & 988 & 988 & 0 \\
3111 & 10 & 610 & 610 & 0 \\
222 & 5 & 298 & 298 & 0 \\
2211 & 9 & 523 & 523 & 0 \\
21111 & 5 & 287 & 287 & 0 \\
111111 & 1 & 56 & 56 & 0 \\
\bottomrule
\end{tabular}
\medskip
\caption{Ranks for BTQ identities in degree 6}
\label{btq6ranks}
\end{table}

%%%%%%%%%%%%%%%%%%%%%%%%%%%%%%%%%%%%%%%%%%%%%%%%%%%%%%%%%%%%%%%%%%%%%%%%%%%%%%%%%%%%

\section{Conclusions and open problems}

In this work we give the defining identities for the tangent algebras to monoassociative
local analytic loops, completing the research of Akivis and many other geometers and
algebraists on nonassociative structures obtained from hexagonal three-webs. However, 
there also appear some new concepts and methods that can be interesting to study into 
more detail. Let us briefly describe some of them.

When defining degree $d$ Sabinin algebras, one natural question arises: do we have examples
of such structures? Which interesting properties do they have? In particular, we are
concerned with the Schreier property of a variety of algebras, that is, if every subalgebra
of a free algebra is also free.
Degree 2 Sabinin algebras are anticommutative algebras, and this is a Schreier variety
by a theorem of Shirshov \cite{Shirshov}.
Degree 3 Sabinin algebras are Akivis algebras, and this variety is also Schreier by a theorem of Shestakov and
Umirbaev \cite{ShestakovUmirbaev}.
Sabinin algebras also satisfy this property by a recent theorem of Chibrikov \cite{Chibrikov}.
It is an interesting open problem to determine whether the
variety of degree $d$ Sabinin algebras also satisfies the Schreier property.

The works of P\'erez-Izquierdo and Shestakov \cite{PIBol,PIS} generalized
the Poincar\'e-Birkhoff-Witt theorem from Lie algebras to Malcev and Bol algebras and
constructed universal enveloping nonassociative algebras for these structures.
We expect that the methods of those papers can be extended to construct universal envelopes
of tangent algebras for monoassociative loops.

As a generalization of the computations we have described in the previous sections, it seems natural
to find simple forms of operations of higher arity analogous to the quaternators in degree 4.
In general, multilinear $n$-ary operations are primitive elements of degree $n$ in the free
nonassociative algebra. We choose them to be the simplest generators for the $S_n$-module of
the multilinear primitive elements in degree $n$ that are not consequences of the operations of
lower arity.
We have found these operations in degree 5.

\begin{definition}
The three \textbf{quinquenators} are these quinquenary operations:
  \begin{alignat*}{2}
  \text{\FiveFlowerOpen}_1(a,b,c,d,e)
  &=
  ( (ab)c, d, e )
  - ( ab, d, e ) c
  - ( ac, d, e ) b
  - ( bc, d, e ) a
  \\
  &\qquad\qquad\qquad
  + ( (a,d,e) c ) b
  + ( (b,d,e) c ) a
  + ( (c,d,e) b ) a,
  \\
  \text{\FiveFlowerOpen}_2(a,b,c,d,e)
  &=
  ( a, (bc)d, e )
  - ( a, bc, e ) d
  - ( a, bd, e ) c
  - ( a, cd, e ) b
  \\
  &\qquad\qquad\qquad
  + ( (a,b,e) d ) c
  + ( (a,c,e) d ) b
  + ( (a,d,e) c ) b,
  \\
  \text{\FiveFlowerOpen}_3(a,b,c,d,e)
  &=
  ( a, b, (cd)e )
  - ( a, b, cd ) e
  - ( a, b, ce ) d
  - ( a, b, de ) c
  \\
  &\qquad\qquad\qquad
  + ( (a,b,c) e ) d
  + ( (a,b,d) e ) c
  + ( (a,b,e) d ) c.
  \end{alignat*}
\end{definition}

The expressions of these elements have an interesting meaning: their vanishing says that we can
reduce a cubic term in the first, second, or third argument of an associator
into a linear combination of terms which have lower degree in that position of the associator.

We are also concerned with the existence of special identities of higher degree in BTQ algebras:
identities which follow from the relations satisfied by the operations defined in the free
power-associative algebra, but which do not follow from the defining identities for BTQ algebras.
We have proved computationally that such identities exist in degree 6, but we do not yet have
an explicit formula.

%%%%%%%%%%%%%%%%%%%%%%%%%%%%%%%%%%%%%%%%%%%%%%%%%%%%%%%%%%%%%%%%%%%%%%%%%%%%%%%%%%%%

\section*{Acknowledgements}

We thank Michael Kinyon for helping us understand the literature on three-webs and analytic loops.
Murray Bremner was supported by a Discovery Grant from NSERC, the Natural
Sciences and Engineering Research Council of Canada. Sara Madariaga was
supported by the Spanish Ministerio de Educaci\'on y Ciencia (FEDER MTM
201-18370-C04-03 and AP2007-01986) and thanks the Department of Mathematics and
Statistics at the University of Saskatchewan for its hospitality during her
visit from September to December 2011.

%%%%%%%%%%%%%%%%%%%%%%%%%%%%%%%%%%%%%%%%%%%%%%%%%%%%%%%%%%%%%%%%%%%%%%%%%%%%%%%%%%%%

\end{document}